\documentstyle[twoside,amstex]{article}
\input amssym.def
\input amssym.tex
\def\bc{\begin{center}}
\def\ec{\end{center}}
\def\no{\noindent}
\def\hang{\hangindent\parindent}
\def\textindent#1{\indent\llap{[#1]\enspace}\ignorespaces}
\def\re{\par\hang\textindent}
\begin{document}
\thispagestyle{empty} \vspace*{3 true cm} \pagestyle{myheadings}
\markboth {\hfill {\sl H. Chen, H. Kose and Y. Kurtulmaz}\hfill}
{\hfill{\sl Strongly clean triangular matrix rings}\hfill}
\vspace*{-1.5 true cm} \bc{\Large\bf Strongly clean triangular
matrix rings\\
\vskip2mm with endomorphisms}\ec

\vskip4mm \bc{\bf H. Chen, H. Kose and Y. Kurtulmaz}\ec

\vskip10mm
\begin{minipage}{120mm}
\bc {\bf Abstract}\ec

\vskip4mm A ring $R$ is strongly clean provided that every element
in $R$ is the sum of an idempotent and a unit that commutate. Let
$T_n(R,\sigma)$ be the skew triangular matrix ring over a local
ring $R$ where $\sigma$ is an endomorphism of $R$. We show that
$T_2(R,\sigma)$ is strongly clean if and only if for any $a\in
1+J(R), b\in J(R)$, $l_a-r_{\sigma(b)}: R\to R$ is surjective.
Further, $T_3(R,\sigma)$ is strongly clean if
$l_{a}-r_{\sigma(b)}, l_{a}-r_{\sigma^2(b)}$ and
$l_{b}-r_{\sigma(a)}$ are surjective for any $a\in U(R),b\in
J(R)$. The necessary condition for $T_3(R,\sigma)$ to be strongly
clean is also obtained. \vskip3mm {\bf 2010 Mathematics subject
classification:}\ \ 16D70, 16E50. \vskip3mm {\bf Keywords and
phrases:}\ \ strongly clean rings, skew triangular matrix rings;
local rings.

\end{minipage}

\vskip20mm \bc{\bf 1. Introduction}\ec

\vskip4mm \no We say that an element $a\in R$ is strongly clean
provided that there exist an idempotent $e\in R$ and a unit $u\in
R$ such that $a=e+u$ and $eu=ue$. A ring $R$ is strongly clean in
case every element in $R$ is strongly clean. Strong cleanness over
commutative rings was extensively studied by many authors from
very different view points (cf. [1-3] and [5-8]). The problem of
deciding the strong cleanness is considerably harder. So far, one
considers strong cleanness only over commutative local rings,
where a ring $R$ is local provided that $R$ has only a maximal
right ideal. As is well known, a ring $R$ is local if and only if
for any $x\in R$, either $x$ or $1-x$ is invertible. In [6], Li
characterizes when $2\times 2$ matrix ring $M_2(R)$ over a
commutative local ring $R$ is strongly clean. The strong cleanness
of triangular matrix rings over such a ring is also investigated
in [1]. For more discussion of strong cleanness , we refer the
reader to [4] and [7].

Let $R$ be a ring, and let $\sigma$ be an endomorphism of $R$. Let
$T_n(R,\sigma)$ be the set of all upper triangular matrices over
the rings $R$. For any $(a_{ij}),(b_{ij})\in T_n(R,\sigma)$, we
define $(a_{ij})+(b_{ij})=(a_{ij}+b_{ij})$, and
$(a_{ij})(b_{ij})=(c_{ij})$ where
$c_{ij}=\sum\limits_{k=i}^{n}a_{ik}\sigma^{k-i}\big(b_{kj}\big)$.
Then $T_n(R,\sigma)$ is a ring under the preceding addition and
multiplication. Clearly, $T_n(R,\sigma)$ will be $T_n(R)$ only
when $\sigma$ is the identity morphism. Let $a\in R$. $l_a: R\to
R$ and $r_a: R\to R$ denote, respectively, the abelian group
endomorphisms given by $l_a(r)=ar$ and $r_a(r)=ra$ for all $r\in
R$. Thus, $l_a-r_b$ is an abelian group endomorphism such that
$(l_a-r_b)(r)=ar-rb$ for any $r\in R$.

The aim of this note is to investigate the strong cleanness over a
noncommutative local ring with an endomorphism. We prove that
$T_2(R,\sigma)$ is strongly clean if and only if for any $a\in
1+J(R), b\in J(R)$, $l_a-r_{\sigma(b)}: R\to R$ is surjective.
Further, $T_3(R,\sigma)$ is strongly clean if
$l_{a}-r_{\sigma(b)}, l_{a}-r_{\sigma^2(b)}$ and
$l_{b}-r_{\sigma(a)}$ are surjective for any $a\in U(R),b\in
J(R)$. The converse is also true if $1_R$ can not be the sum of
two units. These extend the known results of strong cleanness of
matrices over commutative local rings as well.

Throughout, every ring is associative with an identity $1$. $J(R)$
and $U(R)$ will denote, respectively, the Jacobson radical and the
group of units in the ring $R$.

\vskip15mm \bc{\bf 2. The rings $T_2(R,\sigma)$}\ec

\vskip4mm As is well known, the triangular matrix ring $T_2(R)$
over a local ring $R$ is strongly clean if and only if for any
$a\in 1+J(R), b\in J(R)$, $l_a-r_b: R\to R$ is surjective [5,
Theorem 2.2.1]. We extend this result to the skew triangular
matrix ring with an endomorphism.

\vskip4mm \no{\bf Theorem 2.1.}\ \ {\it Let $R$ be a local ring,
and let $\sigma$ be an endomorphism of $R$. Then the following are
equivalent:} \vspace{-.5mm}
\begin{enumerate}
\item [(1)]{\it $T_2(R,\sigma)$ is strongly clean.}
\vspace{-.5mm}
\item [(2)]{\it If $a\in 1+J(R), b\in J(R)$, then $l_a-r_{\sigma(b)}: R\to R$ is surjective.}
\vspace{-.5mm}
\end{enumerate}\no{\it Proof.}\ \ $(1)\Rightarrow (2)$ Let $a\in 1+J(R), b\in
J(R), v\in R$. Then $A=\left(
\begin{array}{cc}
a&-v\\
0&b
\end{array}
\right)\in T_2(R,\sigma)$. By hypothesis, there exists an
idempotent $E=\left(
\begin{array}{cc}
e&x\\
0&f
\end{array}
\right)\in T_2(R,\sigma)$ such that $A-E\in
U\big(T_2(R,\sigma)\big)$ and $AE=EA$. Since $R$ is local, all
idempotents in $R$ are $0$ and $1$. Thus, we see that $e=0,f=1$;
otherwise, $A-E\not\in U\big(T_2(R,\sigma)\big)$. So $E=\left(
\begin{array}{cc}
0&x\\
0&1
\end{array}
\right)$. It follows from $AE=EA$ that $v+x\sigma(b)=ax$, and so
$ax-v=x\sigma(b)$. Therefore we conclude that $l_a-r_{\sigma(b)}:
R\to R$ is surjective.

$(2)\Rightarrow (1)$ Let $A=\left(
\begin{array}{cc}
a&v\\
0&b
\end{array}
\right)\in T_2(R,\sigma)$. If $a,b\in U(R)$, then $A\in
U\big(T_2(R,\sigma)\big)$ is strongly clean. If $a,b\in J(R)$,
then $A-I_2\in U\big(T_2(R,\sigma)\big)$; hence, $A\in
T_2(R,\sigma)$ is strongly clean. Assume that $a\in U(R),b\in
J(R)$. If $a-1\in U(R)$, then $A-I_2\in U\big(T_2(R,\sigma)\big)$;
hence, $A\in T_2(R,\sigma)$ is strongly clean. If $a-1\in J(R)$,
by hypothesis, $l_a-r_{\sigma(b)}: R\to R$ is surjective. Thus,
$ax-x\sigma(b)=-v$ for some $x\in R$. Choose $E= \left(
\begin{array}{cc}
0&x\\
0&1
\end{array}
\right)\in T_2(R,\sigma)$. Then $E=E^2\in T_2(R,\sigma)$. In
addition, $AE=EA$ and $A-E\in U\big(T_2(R,\sigma)\big)$; hence,
$A\in T_2(R,\sigma)$ is strongly clean. Assume that $a\in
J(R),b\in U(R)$. If $b-1\in U(R)$, then $A-I_2\in
U\big(T_2(R,\sigma)\big)$; hence, $A\in T_2(R,\sigma)$ is strongly
clean. If $b-1\in J(R)$, by hypothesis, $l_{1-a}-r_{\sigma(1-b)}:
R\to R$ is surjective. Thus, $(1-a)x-x\sigma(1-b)=v$ for some
$x\in R$. As $\sigma$ is an endomorphism of $R$,
$\sigma(1-b)=\sigma(1)-\sigma(b)=1-\sigma(b)$. Hence,
$ax-x\sigma(b)=-v$. Choose $E= \left(
\begin{array}{cc}
1&x\\
0&0
\end{array}
\right)\in T_2(R,\sigma)$. Then $E=E^2\in T_2(R,\sigma)$. In
addition, $AE=EA$ and $A-E\in U\big(T_2(R,\sigma)\big)$; hence,
$A\in T_2(R,\sigma)$ is strongly clean. Therefore we conclude that
$A\in T_2(R,\sigma)$ is strongly clean in any case.\hfill$\Box$

\vskip4mm Following Diesl [5], a local ring $R$ is bleached
provided that for any $a\in U(R),b\in J(R)$, $l_a-r_b,l_b-r_a$ are
both surjective.

\vskip4mm \no{\bf Corollary 2.2.}\ \ {\it Let $R$ be a local ring,
and let $\sigma$ be an endomorphism of $R$. If $R$ is bleached,
then $T_2(R,\sigma)$ is strongly clean.}\vskip4mm \no{\it Proof.}\
\ Let $a\in 1+J(R), b\in J(R)$. Then $1-a\in J(R), 1-b\in
1+J(R)\subseteq U(R)$. This implies that $\sigma(1-b)\in U(R)$. By
hypothesis, $l_{1-a}-r_{\sigma(1-b)}: R\to R$ is surjective. For
any $v\in R$, we can find some $x\in R$ such that
$(1-a)x-x\sigma(1-b)=-v$. That is, $ax-x\sigma(b)=v$. This implies
that $l_a-r_{\sigma(b)}: R\to R$ is surjective. According to
Theorem 2.1, $T_2(R,\sigma)$ is strongly clean.\hfill$\Box$

\vskip4mm \no{\bf Corollary 2.3.}\ \ {\it Let $R$ be a local ring,
and let $\sigma$ be an endomorphism of $R$. If $J(R)$ is nil, then
$T_2(R,\sigma)$ is strongly clean.}\vskip2mm \no{\it Proof.}\ \
Let $a\in U(R), b\in J(R)$. Then we can find some $n\in {\Bbb N}$
such that $b^n=0$. For any $v\in R$, we choose
$x=\big(l_{a^{-1}}+l_{a^{-2}}r_{b}+ \cdots
+l_{a^{-n}}r_{b^{n-1}}\big)(v)$. One easily checks that
$$\begin{array}{lll}
&&\big(l_a-r_{b}\big)(x)\\
&=&\big(l_a-r_{b}\big)\big(l_{a^{-1}}+l_{a^{-2}}r_{b}+
\cdots +l_{a^{-n}}r_{b^{n-1}}\big)(v)\\
&=&\big(v+a^{-1}vb+\cdots +a^{-n+1}vb^{n-1}\big)-\big(a^{-1}vb+\cdots +a^{-n}vb^n\big) \\
&=&v. \end{array}$$ This shows that $l_a-r_{b}: R\to R$ is
surjective. Likewise, $l_b-r_{a}: R\to R$ is surjective. This
implies that $R$ is bleached. In light of Corollary 2.2,
$T_2(R,\sigma)$ is strongly clean.\hfill$\Box$

\vskip4mm \no{\bf Example 2.4.}\ \ {\it Let ${\Bbb Z}_{p^n}={\Bbb
Z}/p^n{\Bbb Z} (p~\mbox{is prime}, n\in {\Bbb N})$, and let
$\sigma$ be an endomorphism of ${\Bbb Z}_{p^n}$. Then $T_2({\Bbb
Z}_{p^n},\sigma)$ is strongly clean. As ${\Bbb Z}_{p^n}$ is a
local ring with the Jacobson radical $p{\Bbb Z}_{p^n}$. Obviously,
$J\big({\Bbb Z}_{p^n}\big)$ is nil, and we are done by Corollary
2.3.}\hfill$\Box$

\vskip4mm \no{\bf Example 2.5.}\ \ {\it Let ${\Bbb Z}_{4}={\Bbb
Z}/4{\Bbb Z}$, let $R=\{ \left(
\begin{array}{cc}
a&b\\
0&a
\end{array}
\right)~|~a,b\in {\Bbb Z}_{4}\}$, and let $\sigma: R\to R, \left(
\begin{array}{cc}
a&b\\
0&a
\end{array}
\right)\mapsto \left(
\begin{array}{cc}
a&-b\\
0&a
\end{array}
\right)$. Then $T_2(R,\sigma)$ is strongly clean. Obviously,
$\sigma$ is an endomorphism of $R$. It is easy to check that
$J(R)=\{ \left(
\begin{array}{cc}
a&b\\
0&a
\end{array}
\right)~|~a\in 2{\Bbb Z}_{4}, b\in {\Bbb Z}_{4}\}$, and then
$R/J(R)\cong {\Bbb Z}_{2}$ is a field. Thus, $R$ is a local ring.
In addition, $\big(J(R)\big)^4=0$, thus $J(R)$ is nil. Therefore
we are through from Corollary 2.3.}\hfill$\Box$

\vskip4mm Let $\sigma$ be an endomorphism of ${\Bbb
Z}_{3^n}[x]/(x^2+x+1)$. Analogously, $T_2\big({\Bbb
Z}_{3^n}[x]/(x^2+x+1),\sigma\big)$ is strongly clean.

We say that an element $a\in R$ is very clean provided that for
any $x\in R$ there exists an idempotent $e$ such that $ex=xe$ and
either $x-e\in U(R)$ or $x+e\in U(R)$. A ring $R$ is very clean in
case every element in $R$ is very clean. Every clean ring maybe
not strongly clean. For instance, ${\Bbb Z}_{(3)}\bigcap {\Bbb
Z}_{(5)}$ is a very clean ring, but it is not strongly clean.

\vskip4mm \no{\bf Proposition 2.6.}\ \ {\it Let $R$ be a local
ring, and let $\sigma$ be an endomorphism of $R$. Then the
following are equivalent:} \vspace{-.5mm}
\begin{enumerate}
\item [(1)]{\it $T_2(R,\sigma)$ is very clean.}
\vspace{-.5mm}
\item [(2)]{\it $2\in U(R)$ or $T_2(R,\sigma)$ is strongly clean.}
\vspace{-.5mm}
\end{enumerate}\no{\it Proof.}\ \ $(1)\Rightarrow (2)$ Suppose that $2\in J(R)$. Let $a\in 1+J(R), b\in
J(R), v\in R$. Then $A=\left(
\begin{array}{cc}
a&-v\\
0&b
\end{array}
\right)\in T_2(R,\sigma)$. By hypothesis, there exists an
idempotent $E=\left(
\begin{array}{cc}
e&x\\
0&f
\end{array}
\right)\in T_2(R,\sigma)$ such that $A+ E$ or $A-E\in
U\big(T_2(R,\sigma)\big)$ and $AE=EA$. Since $R$ is local, all
idempotents in $R$ are $0$ and $1$.

If $A-E\in U\big(T_2(R,\sigma)\big)$, then we see that $e=0,f=1$;
otherwise, $A-E\not\in U\big(T_2(R,\sigma)\big)$. So $E=\left(
\begin{array}{cc}
0&x\\
0&1
\end{array}
\right)$. It follows from $AE=EA$ that $v+x\sigma(b)=ax$, and so
$ax-v=x\sigma(b)$. Therefore we conclude that $l_a-r_{\sigma(b)}:
R\to R$ is surjective.

If $A+E\in U\big(T_2(R,\sigma)\big)$, then we see that $f=1$;
otherwise, $A-E\not\in U\big(T_2(R,\sigma)\big)$. If $e=0$, then
$E=\left(
\begin{array}{cc}
0&x\\
0&1
\end{array}
\right)$. It follows from $AE=EA$ that $v+x\sigma(b)=ax$, and so
$ax-v=x\sigma(b)$. Therefore we conclude that $l_a-r_{\sigma(b)}:
R\to R$ is surjective. If $e=1$, then $2\in U(R)$, a
contradiction. Therefore $T_2(R,\sigma)$ is strongly clean by
Lemma 2.1.

$(2)\Rightarrow (1)$ If $T_2(R,\sigma)$ is strongly clean, then
$T_2(R,\sigma)$ is very clean. Now we assume that $2\in U(R)$. Let
$A=\left(
\begin{array}{cc}
a&v\\
0&b
\end{array}
\right)\in T_2(R,\sigma)$. If $a,b\in U(R)$, then $A\in
U\big(T_2(R,\sigma)\big)$ is very clean. If $a,b\in J(R)$, then
$A-I_2\in U\big(T_2(R,\sigma)\big)$; hence, $A\in T_2(R,\sigma)$
is very clean. Assume that $a\in U(R),b\in J(R)$. If $a-1\in
U(R)$, then $A-I_2\in U\big(T_2(R,\sigma)\big)$; hence, $A\in
T_2(R,\sigma)$ is very clean.

If $a-1\in J(R)$, then $A+I_2\in U\big(T_2(R,\sigma)\big)$. Hence,
$A$ is very clean.

Assume that $a\in J(R),b\in U(R)$. If $b-1\in U(R)$, then
$A-I_2\in U\big(T_2(R,\sigma)\big)$; hence, $A\in T_2(R,\sigma)$
is very clean.

If $b-1\in J(R)$, then $A+I_2\in U\big(T_2(R,\sigma)\big)$; hence,
$A\in T_2(R,\sigma)$ is very clean. Therefore we conclude that
$A\in T_2(R,\sigma)$ is very clean in any case.\hfill$\Box$

\vskip4mm \no{\bf Example 2.7.}\ \ {\it Let ${\Bbb Z}_{(3)}=\{
\frac{m}{n}~| m,n\in {\Bbb Z}, 3 \nmid n\}$. Then ${\Bbb Z}_{(3)}$
is a local ring in which $2\in U(R)$. In view of Proposition 2.6,
$T_2({\Bbb Z}_{(3)},\sigma)$ is very clean for any endomorphism
$\sigma$ of $R$.}

\vskip15mm \bc{\bf 3. The rings $T_3(R,\sigma)$}\ec

\vskip4mm The goal of this section is to investigate strong
cleanness of $3\times 3$ skew triangular matrix rings with
endomorphisms over a local ring.

\vskip4mm \no{\bf Theorem 3.1.}\ \ {\it Let $R$ be a local ring,
and let $\sigma$ be an endomorphism of $R$. If
$l_{a}-r_{\sigma(b)}, l_{a}-r_{\sigma^2(b)}$ and
$l_{b}-r_{\sigma(a)}$ are surjective for any $a\in U(R),b\in
J(R)$, then $T_3(R,\sigma)$ is strongly clean.}\vskip2mm\no{\it
Proof.}\ \ Let $A=(a_{ij})\in T_3(R,\sigma)$.

Case 1. $a_{11}, a_{22}, a_{33}\in J(R)$. Then $A=I_3+(A-I_3)$,
and so $A-I_3\in U(T_3(R,\sigma))$. Then $A\in T_3(R,\sigma)$ is
strongly clean.

Case 2. $a_{11}\in U(R), a_{22}, a_{33}\in J(R)$. By hypothesis,
we can find some $e_{12}\in R$ such that
$a_{11}e_{12}-e_{12}\sigma(a_{22})=-a_{12}$. Further, we can find
some $e_{13}\in R$ such that
$a_{11}e_{13}-e_{13}\sigma^2(a_{33})=e_{12}\sigma(a_{23})-a_{13}$.
Choose $E=\left(
\begin{array}{ccc}
0&e_{12}&e_{13}\\
0&1&0\\
0&0&1
\end{array}
\right)\in T_3(R,\sigma)$. Then $E=E^2$, and $A=E+(A-E)$, where
$A-E\in U(T_3(R,\sigma))$.In addition,
$$\begin{array}{c}
EA=\left(
\begin{array}{ccc}
0&e_{12}\sigma(a_{22})&e_{12}\sigma(a_{23})+e_{13}\sigma^2(a_{33})\\
0&a_{22}&a_{23}\\
0&0&a_{33}
\end{array}
\right),\\AE=\left(
\begin{array}{ccc}
0&a_{11}e_{12}+a_{12}&a_{11}e_{13}+a_{13}\\
0&a_{22}&a_{23}\\
0&0&a_{33}\end{array} \right), \end{array}$$ and so $EA=AE$. Hence
$A\in T_3(R,\sigma)$ is strongly clean.

Case 3. $a_{11}\in J(R), a_{22}\in U(R), a_{33}\in J(R)$. Clearly
$\sigma(a_{22})\in U(R).$ By hypothesis, we can find some
$e_{12}\in R$ such that
$a_{11}e_{12}-e_{12}\sigma(a_{22})=a_{12}$. Further, we have some
$e_{23}\in R$ such that
$a_{22}e_{23}-e_{23}\sigma(a_{33})=-a_{23}$. It follows that
$-a_{11}e_{12}\sigma(e_{23})+a_{12}\sigma(e_{23})=-e_{12}\sigma(a_{22})\sigma(e_{23})=e_{12}\sigma(a_{23})-e_{12}\sigma(e_{23})\sigma^2(a_{33}).$
Choose $E=\left(
\begin{array}{ccc}
1&e_{12}&-e_{12}\sigma(e_{23})\\
0&0&e_{23}\\
0&0&1
\end{array}
\right)\in T_3(R,\sigma)$. Then $E=E^2$, and $A=E+(A-E)$, where
$A-E\in U(T_3(R,\sigma))$. In addition,
$$\begin{array}{c}
EA=\left(
\begin{array}{ccc}
a_{11}&a_{12}+e_{12}\sigma(a_{22})&a_{13}+e_{12}\sigma(a_{23})-e_{12}\sigma(e_{23})\sigma^2(a_{33})\\
0&0&e_{23}\sigma(a_{33})\\
0&0&a_{33}
\end{array}
\right),\\
AE=\left(
\begin{array}{ccc}
a_{11}&a_{11}e_{12}&-a_{11}e_{12}\sigma(e_{23})+a_{12}\sigma(e_{23})+a_{13}\\
0&0&a_{22}e_{23}+a_{23}\\
0&0&a_{33}
\end{array}
\right), \end{array}$$ and so $EA=AE$. Hence $A\in T_3(R,\sigma)$
is strongly clean.

Case 4. $a_{11}, a_{22}\in J(R), a_{33}\in U(R)$. By hypothesis,
we can find some $e_{23}\in R$ such that
$a_{22}e_{23}-e_{23}\sigma(a_{33})=a_{23}$. Clearly,
$\sigma(a_{33})\in U(R)$. Thus, we can find some $e_{13}\in R$
such that
$a_{11}e_{13}-e_{13}\sigma^2(a_{33})=a_{13}-a_{12}\sigma(e_{23})$.
Choose $E=\left(
\begin{array}{ccc}
1&0&e_{13}\\
0&1&e_{23}\\
0&0&0
\end{array}
\right)\in T_3(R,\sigma)$. Then $E=E^2$, and $A=E+(A-E)$, where
$A-E\in U(T_3(R,\sigma))$. In addition,
$$\begin{array}{c}
EA=\left(
\begin{array}{ccc}
a_{11}&a_{12}&a_{13}+e_{13}\sigma^2(a_{33})\\
0&a_{22}&a_{23}+e_{23}\sigma(a_{33})\\
0&0&0
\end{array}
\right),\\
AE=\left(
\begin{array}{ccc}
a_{11}&a_{12}&a_{11}e_{13}+a_{12}\sigma(e_{23})\\
0&a_{22}&a_{22}e_{23}\\
0&0&0
\end{array}
\right), \end{array}$$ and so $EA=AE$. Hence $A\in T_3(R,\sigma)$
is strongly clean.

Case 5. $a_{11}\in J(R), a_{22}, a_{33}\in U(R)$. By hypothesis,
we can find some $e_{12}\in R$ such that
$a_{11}e_{12}-e_{12}\sigma(a_{22})=a_{12}$. Further, we can find
some $e_{13}\in R$ such that
$a_{11}e_{13}-e_{13}\sigma^2(a_{33})=a_{13}+e_{12}\sigma(e_{23})$.
Choose $E=\left(
\begin{array}{ccc}
1&e_{12}&e_{13}\\
0&0&0\\
0&0&0
\end{array}
\right)\in T_3(R,\sigma)$. Then $E=E^2$, and $A=E+(A-E)$, where
$A-E\in U(T_3(R,\sigma))$. In addition,
$$\begin{array}{c}
EA=\left(
\begin{array}{ccc}
a_{11}&a_{12}+e_{12}\sigma(a_{22})&a_{13}+e_{12}\sigma(a_{23})+e_{13}\sigma^2(a_{33})\\
0&0&0\\
0&0&0
\end{array}
\right),\\
AE=\left(
\begin{array}{ccc}
a_{11}&a_{11}e_{12}&a_{11}e_{13}\\
0&0&0\\
0&0&0
\end{array}
\right), \end{array}$$ and so $EA=AE$. Hence $A\in T_3(R,\sigma)$
is strongly clean.

Case 6. $a_{11}\in U(R), a_{22}\in J(R), a_{33}\in U(R)$. By
hypothesis, we can find some $e_{23}\in R$ such that
$a_{22}e_{23}-e_{23}\sigma(a_{33})=a_{23}$. Further, we can find
some $e_{12}\in R$ such that
$a_{11}e_{12}-e_{12}\sigma(a_{22})=-a_{12}$. It is easy to verify
that
$$e_{12}\sigma(a_{23})+e_{12}\sigma(e_{23})\sigma^2(a_{33})=e_{12}\sigma(a_{22}e_{23})=a_{11}e_{12}\sigma(e_{23})+a_{12}\sigma(e_{23}).$$
Choose $E=\left(
\begin{array}{ccc}
0&e_{12}&e_{12}\sigma(e_{23})\\
0&1&e_{23}\\
0&0&0
\end{array}
\right)\in T_3(R,\sigma)$. Then $E=E^2$, and $A=E+(A-E)$, where
$A-E\in U(T_3(R,\sigma))$. In addition, $$\begin{array}{c}
EA=\left(
\begin{array}{ccc}
0&e_{12}\sigma(a_{22})&e_{12}\sigma(a_{23})+e_{12}\sigma(e_{23})\sigma^2(a_{33})\\
0&a_{22}&a_{23}+e_{23}\sigma(a_{33})\\
0&0&0
\end{array}
\right),\\
AE=\left(
\begin{array}{ccc}
0&a_{11}e_{12}+a_{12}&a_{11}e_{12}\sigma(e_{23})+a_{12}\sigma(e_{23})\\
0&a_{22}&a_{22}e_{23}\\
0&0&0
\end{array}
\right) \end{array}$$ and so $EA=AE$. Hence $A\in T_3(R,\sigma)$
is strongly clean.

Case 7. $a_{11}, a_{22}\in U(R), a_{33}\in J(R)$. By hypothesis,
we can find some $e_{23}\in R$ such that
$a_{22}e_{23}-e_{23}\sigma(a_{33})=-a_{23}$. Further, we can find
some $e_{13}\in R$ such that
$a_{11}e_{13}-e_{13}\sigma^2(a_{33})=-a_{13}-a_{12}\sigma(e_{23})$.
Choose $E=\left(
\begin{array}{ccc}
0&0&e_{13}\\
0&0&e_{23}\\
0&0&1
\end{array}
\right)\in T_3(R,\sigma)$. Then $E=E^2$, and $A=E+(A-E)$, where
$A-E\in U(T_3(R,\sigma))$. In addition, $$\begin{array}{c}
EA=\left(
\begin{array}{ccc}
0&0&e_{13}\sigma^2(a_{33})\\
0&0&e_{23}\sigma(a_{33})\\
0&0&a_{33}
\end{array}
\right),\\
AE=\left(
\begin{array}{ccc}
0&0&a_{11}e_{13}+a_{12}\sigma(e_{23})+a_{13}\\
0&0&a_{22}e_{23}+a_{23}\\
0&0&a_{33}
\end{array}
\right), \end{array}$$ and so $EA=AE$. Hence $A\in T_3(R,\sigma)$
is strongly clean.

Case 8. $a_{11}, a_{22}, a_{33}\in U(R)$. Then $A=0+A$, where
$A\in U(T_3(R,\sigma))$. Hence $A\in T_3(R,\sigma)$ is strongly
clean.

Therefore we conclude that $T_3(R,\sigma)$ is strongly
clean.\hfill$\Box$

\vskip4mm \no{\bf Corollary 3.2.}\ \ {\it Let $R$ be a local ring,
and let $\sigma$ be an endomorphism of $R$. If $J(R)$ is nil, then
$T_3(R,\sigma)$ is strongly clean.}\vskip2mm \no{\it Proof.}\ \
Let $a\in U(R), b\in J(R)$. Then we can find some $n\in {\Bbb N}$
such that $b^n=0$; hence, $\big(\sigma(b)\big)^n=0$. For any $v\in
R$, we choose $x=\big(l_{a^{-1}}+l_{a^{-2}}r_{\sigma(b)}+ \cdots
+l_{a^{-n}}r_{(\sigma(b))^{n-1}}\big)(v)$. One easily checks that
$$\begin{array}{lll}
&&\big(l_a-r_{\sigma(b)}\big)(x)\\
&=&\big(l_a-r_{\sigma(b)}\big)\big(l_{a^{-1}}+l_{a^{-2}}r_{\sigma(b)}+
\cdots +l_{a^{-n}}r_{\big(\sigma(b)\big)^{n-1}}\big)(v)\\
&=&\big(v+a^{-1}v\sigma(b)+\cdots +a^{-n+1}v\big(\sigma(b)\big)^{n-1}\big)-\big(a^{-1}v\sigma(b)+\cdots +a^{-n}v\big(\sigma(b)\big)^n\big) \\
&=&v. \end{array}$$ This shows that $l_a-r_{\sigma(b)}: R\to R$ is
surjective. Likewise, $l_a-r_{\sigma^2(b)}, l_b-r_{\sigma(a)}:
R\to R$ are surjective. Therefore $T_3(R,\sigma)$ is strongly
clean by Theorem 3.1.\hfill$\Box$

\vskip4mm \no{\bf Corollary 3.3.}\ \ {\it Let $R$ be a bleached
ring, and let $\sigma$ be an endomorphism of $R$. If
$\sigma\big(J(R)\big)\subseteq J(R)$, then $T_3(R,\sigma)$ is
strongly clean.}\vskip2mm \no{\it Proof.}\ \ Let $a\in U(R), b\in
J(R)$. Since $\sigma\big(J(R)\big)\subseteq J(R)$, we see that
$\sigma(b), \sigma^2(b)\in J(R)$. Clearly, $\sigma(a)\in U(R)$. As
$R$ is bleached, we get $l_{a}-r_{\sigma(b)},
l_{a}-r_{\sigma^2(b)}$ and $l_{b}-r_{\sigma(a)}$ are surjective .
According to Theorem 3.1, the result follows.\hfill$\Box$

\vskip4mm \no{\bf Example 3.4.}\ \ {\it Let $R=\big({\Bbb
Z}[x]\big)_{(x)}$ be the localization of ${\Bbb Z}[x]$ at the
prime ideal $(x)$. Define $\sigma: R\to R$ given by $\sigma\big(
\frac{f(x)}{g(x)}\big)=\frac{f(0)}{g(0)}$. Then $\sigma: R\to R$
is an endomorphism over $R$. Moreover, $R$ is a bleached local
ring and $\sigma\big(J(R)\big)\subseteq J(R)$. In light of
Corollary 3.3, $T_3(R,\sigma)$ is strongly clean.}\hfill$\Box$

\vskip10mm \bc{\bf 4. A necessary condition}\ec

\vskip4mm This section is concerned with the necessary condition
on a local ring $R$ under which the skew triangular matrix ring
$T_3(R,\sigma)$ is strongly clean.

\vskip4mm \no{\bf Theorem 4.1.}\ \ {\it Let $R$ be a local ring,
and let $\sigma$ be an endomorphism of $R$. If $T_3(R,\sigma)$ is
strongly clean, then $l_{a}-r_{\sigma(b)}, l_{a}-r_{\sigma^2(b)}$
and $l_{b}-r_{\sigma(a)}$ are surjective for any $a\in 1+J(R),b\in
J(R)$.}\vskip2mm \no{\it Proof.}\ \ Let $a\in 1+J(R),b\in J(R)$.
Choose $E=\left(
\begin{array}{cc}
I_2&\\
&0\\
\end{array}
\right)$. Then $T_2(R,\sigma)\cong ET_3(R,\sigma)E$, and so
$T_2(R,\sigma)$ is strongly clean. According to Theorem 2.1,
$l_{a}-r_{\sigma(b)}$ is surjective. As $1-b\in 1+J(R)$ and
$1-a\in J(R)$, by the preceding discussion, we see that
$l_{1-b}-r_{\sigma(1-a)}: R\to R$ is surjective. For any $v\in R$,
we can find some $x\in R$ such that $(1-b)x-x\sigma(1-a)=-v$. This
implies that $bx-x\sigma(a)=v$. Hence, $l_b-r_{\sigma(a)}: R\to R$
is surjective.

Let $v\in R$ and let $A=\left(
\begin{array}{ccc}
b&0&v\\
0&b&0\\
0&0&a \end{array} \right)\in T_3(R,\sigma)$. Then we can find an
idempotent $E=(e_{ij})\in T_3(R,\sigma)$ such that $A-E\in U\big(
T_3(R,\sigma)\big)$ and $EA=AE$. This implies that
$e_{11},e_{22},e_{33}\in R$ are all idempotents. As $a\in 1+J(R),
b\in J(R)$, we get $e_{11}=1, e_{22}=1$ and $e_{33}=0$; otherwise,
$A-E$ is not invertible. As $E=E^2$, we see that $E=\left(
\begin{array}{ccc}
1&0&e_{13}\\
0&1&e_{23}\\
0&0&0 \end{array} \right)$ for some $e_{13}, e_{23}\in R$.
Moreover, we see that
$$\left(
\begin{array}{ccc}
b&0&be_{13}\\
0&b&be_{23}\\
0&0&0 \end{array} \right)=AE=EA=\left(
\begin{array}{ccc}
b&0&v+e_{13}\sigma^2(a)\\
0&b&e_{23}\sigma(a)\\
0&0&0 \end{array} \right),$$ and so $be_{13}-e_{13}\sigma^2(a)=v$.
This means that $l_b-r\sigma^2(a): R\to R$ is surjective. As
$1-a\in J(R)$ and $1-b\in 1+J(R)$, by the preceding discussion,
$l_{1-a}-r_{\sigma^2(1-b)}: R\to R$ is surjective. Thus, we can
find some $x\in R$ such that $(1-a)x-x\sigma^2(1-b)=-v$. This
implies that $ax-x\sigma^2(b)=v$, and so $l_{a}-r_{\sigma^2(b)}$
is surjective, as desired. \hfill$\Box$

\vskip4mm \no{\bf Corollary 4.2.}\ \ {\it Let $R$ be a local ring
in which $1_R$ is not the sum of two units, and let $\sigma$ be an
endomorphism of $R$. Then the following are equivalent:}
\vspace{-.5mm}
\begin{enumerate}
\item [(1)]{\it $T_3(R,\sigma)$ is strongly clean.}
\vspace{-.5mm}
\item [(2)]{\it $l_{a}-r_{\sigma(b)}, l_{a}-r_{\sigma^2(b)}$ are surjective for any $a\in 1+J(R),b\in
J(R)$.} \vspace{-.5mm}
\end{enumerate}\no{\it Proof.}\ \ $(1)\Rightarrow (2)$ is obvious
from Theorem 4.1.

$(2)\Rightarrow (1)$ For any $a\in 1+J(R),b\in J(R)$, as in the
proof of Theorem 4.1, we see that $l_b-r_{\sigma(a)}: R\to R$ is
surjective. Obviously, $1+J(R)=U(R)$. For any $u\in U(R)$, we know
that either $u-1\in J(R)$ or $u-1\in U(R)$. Thus $u\in 1+J(R)$;
otherwise, $1=u+(1-u)$ is the sum of two units, a contradiction.
Therefore $U(R)=1+J(R)$. According to Theorem 3.1, $T_3(R,\sigma)$
is strongly clean. \hfill$\Box$

\vskip4mm \no{\bf Corollary 4.3.}\ \ {\it Let $R$ be a local ring
in which $1$ is not the sum of two units, and let
$\sigma=\sigma^2$ be an endomorphism of $R$.Then the following are
equivalent:} \vspace{-.5mm}
\begin{enumerate}
\item [(1)]{\it $T_2(R,\sigma)$ is strongly clean.}
\vspace{-.5mm}
\item [(2)]{\it $T_3(R,\sigma)$ is strongly clean.}
\vspace{-.5mm}
\item [(3)]{\it $l_{a}-r_{\sigma(b)}$ is surjective for any $a\in 1+J(R),b\in
J(R)$.} \vspace{-.5mm}
\end{enumerate}\no{\it Proof.}\ \ $(1)\Leftrightarrow (3)$ is proved by Theorem 2.1.

$(2)\Leftrightarrow (3)$ is obvious from Corollary
4.2.\hfill$\Box$

\vskip4mm \no{\bf Example 4.4.}\ \ {\it Let ${\Bbb Z}_4={\Bbb
Z}/4{\Bbb Z}$, and let $G=\{1, g\}$ be the abelian group of order
$2$. Let $\sigma: {\Bbb Z}_4G\to {\Bbb Z}_4G, a+bg\to a+b$ for any
$a+bg\in {\Bbb Z}_4G$. Then $T_2({\Bbb Z}_4G,\sigma)$ and
$T_3({\Bbb Z}_4G,\sigma)$ are strongly clean.}\vskip2mm \no{\it
Proof.}\ \ Clearly, ${\Bbb Z}_4$ is a local ring with the Jacobson
radical $2{\Bbb Z}_4=\{ \overline{0},\overline{2}\}$. It is easy
to verify that $a+bg\in U\big({\Bbb Z}_4G\big)$ if and only if
$a+b\in U({\Bbb Z}_4)$. Thus, $J({\Bbb Z}_4G)=\{ a+bg~|~a+b\in
J({\Bbb Z}_4)\}$. If $(a+bg)+(c+dg)=1$, then $a+c=1$ and $b+d=0$,
and so $a+b+c+d=1$. If $a+b and c+d\not\in U({\Bbb Z}_4)$, then
$a+b,c+d=0,2$. Hence, $a+b+c+d=0,2$, a contradiction. Thus,
$a+b\in U({\Bbb Z}_4)$ or $c+d\in U({\Bbb Z}_4)$. That is,
$a+bg\in U\big({\Bbb Z}_4G\big)$ or $c+dg\in U\big({\Bbb
Z}_4G\big)$. This implies that ${\Bbb Z}_4G$ is a local ring. If
$1_{RG}=(a+bg)+(c+dg)$ where $a+bg,c+dg\in U\big({\Bbb
Z}_4G\big)$, then $a+c=1$ and $b+d=0$. This yields that
$a+c+b+d=1$. On the other hand, $a+b,c+d=1,3$; hence,
$a+b+c+d=0,2$, a contradiction. This implies that $1_{{\Bbb
Z}_4G}$ is not the sum of two units. Obviously, $\sigma=\sigma^2$.
As ${\Bbb Z}_4G$ is commutative, the preceding condition $(3)$
holds. According to Corollary 4.3, we are done.\hfill$\Box$

\vskip4mm We end this note by asking a problem: How to
characterize a strongly clean $n$ by $n$ triangular matrix ring
$T_n(R,\sigma)$ for $n\geq 4$?

\vskip15mm \bc{\large \bf References}\ec \vskip 4mm \small{\re{1}
M.S. Ahn, Weakly Clean Rings and Almost Clean rings, Ph.D. Thesis,
The University of Iowa, 2003.

\re{2} H. Chen, On uniquely clean rings,
 {\it Comm. Algebra}, {\bf 39}(2011), 189--198.

\re{3} H. Chen, {\it Rings Related Stable Range Conditions}, {\it
Series in Algebra} {\bf 11}, Hackensack, NJ: World Scientific,
2011.

\re{4} H. Chen; O. Gurgun and H. Kose, Strongly clean matrices
over commutative local rings, {\it J. Algebra Appl.}, 12, No. 1,
Paper No. 1250126, 13 p. (2013).

\re{5} A.J. Diesl, {\it Classes of Strongly Clean Rings}, Ph.D.
Thesis, University of California, Berkeley, 2006.

\re{6} Y. Li, On strongly clean matrix rings, {\it J. Algebra},
{\bf 312}(2007), 397--404.

\re{7} W.K. Nicholson, Clean rings: a survey, {\it Advances in
Ring Theory}, World Sci. Publ., Hackensack, NJ, 2005, 181--198.

\re{8} W.K. Nicholson and Y. Zhou, Rings in which elements are
uniquely the sum of an idempotent and a unit, {\it Glasg. Math.
J.}, {\bf 46}(2004), 227--236.

\vskip10mm
\no H. Chen\\
Department of Mathematics\\
Hangzhou Normal University\\
Hangzhou 310034, China\\
E-mail: huanyinchen@@yahoo.cn

\vskip4mm\no H. Kose\\
Department of Mathematics\\
Ahi Evran University, Kirsehir, Turkey\\
E-mail: handankose@@gmail.com

\vskip4mm\no Y. Kurtulmaz\\
Department of Mathematics\\
Bilkent University, Ankara, Turkey\\
E-mail: yosum@@fen.bilkent.edu.tr

\end{document}